\documentclass[12pt]{article}

\usepackage{amsfonts}
\usepackage{amsthm}
\usepackage[centertags]{amsmath}
\usepackage{tikz-cd}
\usepackage{enumerate}
\usepackage[pdftex]{pict2e}
\usepackage{xcolor}

\usepackage{graphicx}
\usepackage{import}
\usepackage{epstopdf}

\newtheorem{theorem}{Theorem}
\newtheorem*{theorem*}{Theorem}
\newtheorem{lemma}[theorem]{Lemma}

\theoremstyle{plain}

\begin{document}

\title{The shortest way to
the geodesics of spheres}

\author{
Mauro Patr\~{a}o\footnote{Departament of Mathematics,
University of Brasília, Brazil. \textit{mpatrao@mat.unb.br }.
}
}

\maketitle

\begin{abstract}
In this paper, we prove, using only elementary geometric arguments and only assuming that the curves are continuous, that the geodesics on a sphere are the minor arcs of the great circles. Our result are valid for any sphere in any inner product space.
\end{abstract}
  


The question of what is the easiest way to obtain the geodesics of spheres is an old one. In \cite{chorlton}, a very short way is obtained for smooth curves, avoiding the use of the calculus of variations. In the present paper, we only assume continuity to obtain not only the (local) uniqueness result, but also to provide the existence result for geodesics of spheres, which is missing in \cite{banks}. Our approach is more direct than of \cite{banks}, since we rely only on rotations and do not use any property of intersection of curves which are valid only for surfaces, and thus our result applies to any sphere in any inner product space $V$.

The length of a given continuous curve $C: [a,b] \to V$ is given by
\[
\ell(C) = \sup \sum_{i=0}^n d(C(t_i),C(t_{i+1}))
\]
where $d(p,q)$ is the distance between the points $p$ and $q$ in $V$ induced by the inner product and the supremum is taken over all finite partitions of $[a,b]$, given by $a = t_0 < t_1 < \cdots < t_n < t_{n+1} = b$. It is immediate from the definition that $\ell(C) \geq d(C(a),C(b))$. Given $a \leq s \leq t \leq b$, we denote by $C_{st}$ the restriction of $C$ to $[s,t]$ and, without loss of generality, we assume that $\ell(C_{st}) > 0$ if $s < t$. Given $C_1: [a_1,b_1] \to V$ and $C_2: [a_2,b_2] \to V$ two continuous curves such that $C_1(b_1) = C_2(a_2)$, we can connect these two curves, defining the following continuous curve $C_1 \cup C_2: [a_1,b_2-a_2+b_1] \to V$, given by
\[
\left\{
\begin{array}{lcll}
(C_1 \cup C_2)(t) & = & C_1(t),& t \in [a_1,b_1] \\
(C_1 \cup C_2)(t) & = & C_2(t-b_1+a_2),& t \in [b_1,b_2-a_2+b_1]
\end{array}
\right.
\]
It is not difficult to prove that
\[
\ell(C_1 \cup C_2) = \ell(C_1) + \ell(C_2)
\]

Denote by $\pi_p q$ the perpendicular projection of a point $q \in S$ on the axis defined by $p$, which is defined by the straight line containing $p$ and the center. The set of points $q \in S$ such that $\pi_p q$ is constant is given by the intersection of $S$ with a plane perpendicular to the axis defined by $p$ and it is a circle in this plane. For any point in this circle, there exists a rotation that takes this point to $q$ while keeping $p$ fixed. When $\pi_p q$ is the center, the circle is called a great circle.

\begin{center}
\begin{picture}(220,180)
\put(-80,0){\includegraphics[scale=.25]{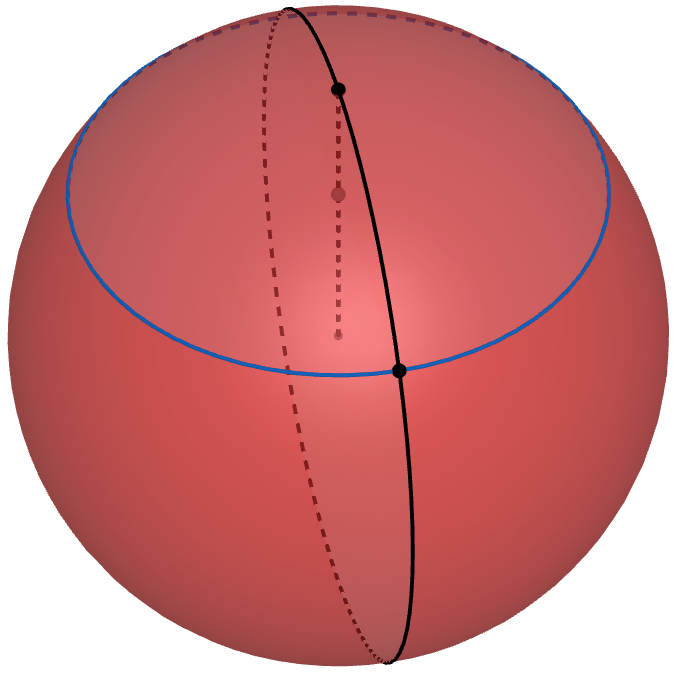}}
\put(120,0){\includegraphics[scale=.29]{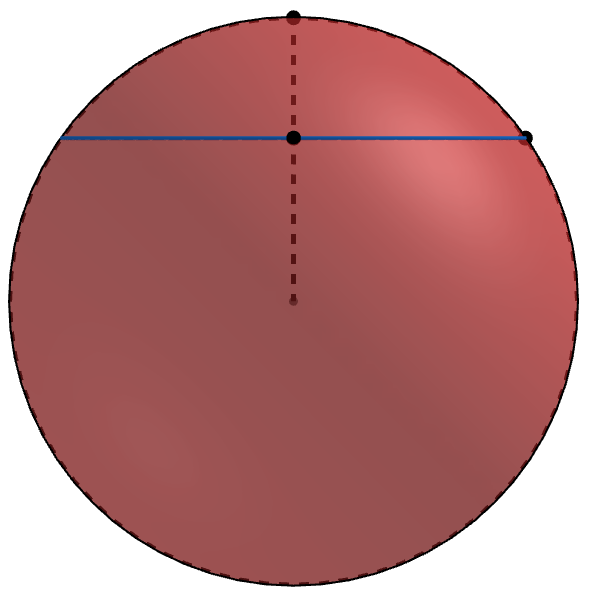}}
\put(10,150){$p$}
\put(25,85){$q$}
\put(-15,120){$\pi_p q$}
\put(200,175){$p$}
\put(280,130){$q$}
\put(185,122){$\pi_p q$}
\end{picture}

\end{center}

It is always possible to put the point $p$ in north pole and to the point $q$ in such way that the points $p$, $q$, and the center are in a plane parallel to the page. The perpendicular projection of $S$ onto the page is a disc whose boundary is a great circle containing $p$ and $q$. The perpendicular projections of the circles such that $\pi_p$ and $\pi_q$ onto the page are, respectively, constant are chords perpendicular to the axis defined by, respectively, $p$ anda $q$. Observe that in order to obtain the perpendicular projection of a point in $S$ onto the axis defined by $p$ or $q$ we can first obtain the perpendicular projection of this point onto the page and then obtain the perpendicular projection of this projection onto the axis defined by $p$ or $q$.

When $C: [a,b] \to V$ is such that $C(a) = p$ and that $C(b) = q$, we say that $C$ is a curve from $p$ to $q$, and, when $C([a,b]) \subset S$, we say that $C$ lies on the sphere $S$ and denote this by $C \subset S$. A curve $G: [a,b] \to V$ from $p$ to $q$ such that $G \subset S$ is called a geodesic between $p$ and $q$ when $G$ has the shortest path in $S$ between $p$ and $q$, i.e., $\ell(G) \leq \ell(C)$, for any curve $C \subset S$ from $p$ to $q$. From now on, we assume that $C$ is a curve from $p$ to $q$ such that $C \subset S$.

\begin{lemma}
If there exist $a < s < t < b$ such that $\pi_p C(s) = \pi_p C(t)$ or that $\pi_q C(s) = \pi_q C(t)$, then $C$ is not a geodesic between $p$ and $q$.
\end{lemma}

\begin{proof}
we can write $C = C_s \cup C_{st} \cup C_t$, where $C_s = C_{as}$ and $C_t = C_{tb}$. If $\pi_p C(s) = \pi_p C(t)$, we can rotate the curve $C_s$ through the axis defined by $p$ obtaining a curve $\widehat{C}_s$ from $p$ to $C(t)$. Now we have that $\widehat{C} = \widehat{C}_s \cup C_t$ is a curve from $p$ to $q$ such that $\widehat{C} \subset S$ and that
\[
\ell(\widehat{C}) = \ell(\widehat{C}_s) + \ell(C_t) < \ell(C_s) + \ell(C_{st}) + \ell(C_t)  = \ell(C)
\]
showing that $C$ is not a geodesic from $p$ to $q$. If $\pi_q C(s) = \pi_q C(t)$, we can rotate the curve $C_t$ through the axis defined by $q$ obtaining a curve $\widetilde{C}_t$ from $C(s)$ to $q$. Now we have that $\widetilde{C} = C_s \cup \widetilde{C}_t$ is a curve from $p$ to $q$ such that $\widetilde{C} \subset S$ and that
\[
\ell(\widetilde{C}) = \ell(C_s) + \ell(\widetilde{C}_t) < \ell(C_s) + \ell(C_{st}) + \ell(C_t)  = \ell(C)
\]
showing again that $C$ is not a geodesic from $p$ to $q$.
\end{proof}

\begin{center}
\begin{picture}(220,180)
\put(-80,0){\includegraphics[scale=.25]{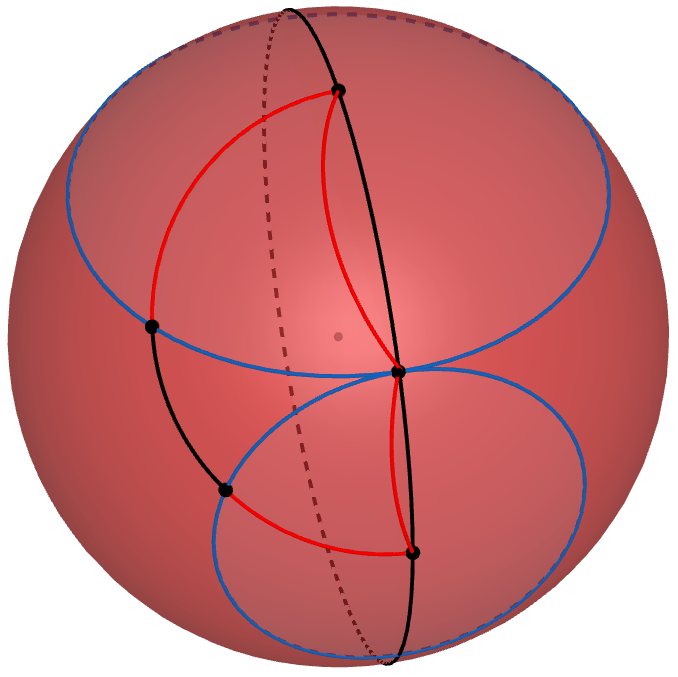}}
\put(120,0){\includegraphics[scale=.29]{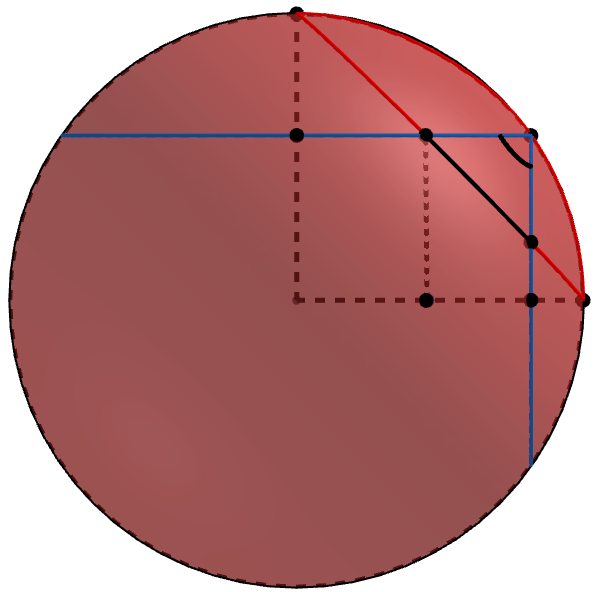}}
\put(10,150){$p$}
\put(25,85){$g$}
\put(30,30){$q$}
\put(-70,80){$C(s)$}
\put(-50,40){$C(t)$}
\put(-45,125){$\textcolor{red}{C_{s}}$}
\put(-10,105){$\textcolor{red}{\widehat{C}_{s}}$}
\put(-53,60){$C_{st}$}
\put(-10,23){$\textcolor{red}{C_{t}}$}
\put(5,45){$\textcolor{red}{\widehat{C}_{t}}$}
\put(200,175){$p$}
\put(280,130){$g$}
\put(295,85){$q$}
\put(262,120){$\alpha$}
\put(219,121){$C(s)$}
\put(185,122){$\pi_p g$}
\put(252,92){$C(t)$}
\put(256,74){$\pi_q g$}
\put(210,73){$\pi_q C(s)$}
\end{picture}

\end{center}

\begin{theorem}
If $C$ is not the minor arc from $p$ to $q$ of a great circle, then $C$ is not a geodesic from $p$ to $q$.
\end{theorem}

\begin{proof}
By the lemma, we may assume that $d(\pi_p C(t),p)$ and $d(\pi_q C(t),q)$ are, respectively, increasing and decreasing functions of $t$. We first consider the case where $p$ and $q$ are not antipodal points, so that there exist a unique great circle containing $p$ and $q$. Let $s$ be such that $C(s)$ is not in this great circle. Since the image by $\pi_p$ of the minor arc from $p$ to $q$ is the segment from $p$ to $\pi_p q$ there exists a unique point $g$ in the great circle such that $\pi_p C(s) = \pi_p g$. Since the image by $\pi_q$ of the minor arc from $p$ to $q$ is the segment from $q$ to $\pi_q p$ there exists some $t$ such that $\pi_q C(t) = \pi_q g$. The angle $\alpha$ between the segment from $C(s)$ to $g$ and the segment from $g$ to $C(t)$ is the supplement of the angle between the axes defined by, respectively, $p$ and $q$, which is positive, since $p$ and $q$ are not antipodal points. Hence $d(\pi_q C(s),q) > d(\pi_q C(t),q)$, which implies that $s < t$. Now we can write $C = C_s \cup C_{st} \cup C_t$, where $C_s = C_{as}$ and $C_t = C_{tb}$. We can rotate the curve $C_s$ through the axis defined by $p$ obtaining a curve $\widehat{C}_s$ from $p$ to $g$ and we can rotate the curve $C_t$ through the axis defined by $q$ obtaining a curve $\widehat{C}_t$ from $g$ to $q$. Now we have that $\widehat{C} = \widehat{C}_s \cup \widehat{C}_t$ is a curve from $p$ to $q$ such that $\widehat{C} \subset S$ and that
\[
\ell(\widehat{C}) = \ell(\widehat{C}_s) + \ell(\widehat{C}_t) < \ell(C_s) + \ell(C_{st}) + \ell(C_t)  = \ell(C)
\]
showing that $C$ is not a geodesic from $p$ to $q$. If $p$ and $q$ are antipodal points, we can consider an intermediate point $C(t)$ such $p$ and $C(t)$ are not antipodal points as well $C(t)$ and $q$. Now we can apply the previous result to $C_{at}$ and to $C_{tb}$.
\end{proof}

\begin{center}
\begin{picture}(230,180)
\put(-80,0){\includegraphics[scale=.25]{sphere-4.png}}
\put(120,0){\includegraphics[scale=.29]{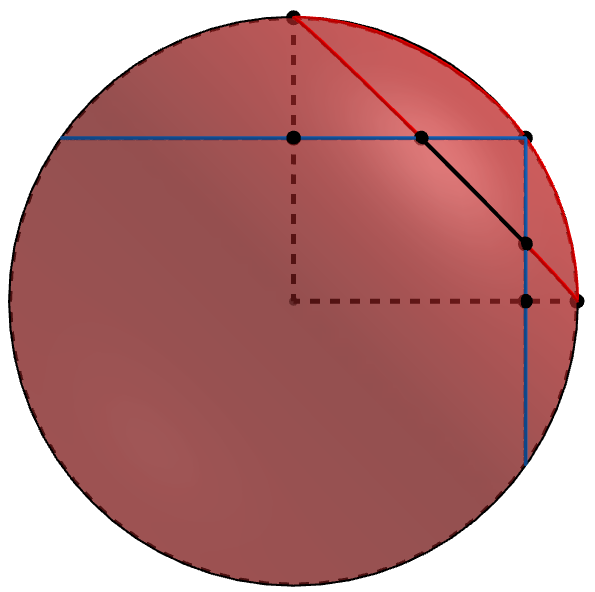}}
\put(10,150){$p$}
\put(20,95){$G(t_{n+1})$}
\put(30,30){$q$}
\put(-70,80){$C(s)$}
\put(-50,40){$C(t)$}
\put(-45,125){$\textcolor{red}{C_{s}}$}
\put(-10,105){$\textcolor{red}{\widehat{C}_{s}}$}
\put(-53,60){$C_{st}$}
\put(-10,23){$\textcolor{red}{C_{t}}$}
\put(5,45){$\textcolor{red}{\widehat{C}_{t}}$}
\put(200,175){$p$}
\put(278,130){$G(t_{n+1})$}
\put(295,85){$q$}
\put(219,121){$C(s)$}
\put(153,121){$\pi_p G(t_{n+1})$}
\put(249,92){$C(t)$}
\put(221,72){$\pi_q G(t_{n+1})$}
\end{picture}

\end{center}

\begin{theorem}
If $G$ is the minor arc from $p$ to $q$ of a great circle, then $G$ is a geodesic from $p$ to $q$.
\end{theorem}

\begin{proof}
First, for any $p,q \in S$ and any curve $C \subset S$ from $p$ to $q$, we prove by induction in $n$ that
\[
\ell(C) \geq \sum_{i=0}^n d(G(t_i),G(t_{i+1}))
\]
where $a = t_0 < t_1 < \cdots < t_n < t_{n+1} = b$ is any finite partition of $[a,b]$. For $n = 0$, we have that
\[
\ell(C) \geq d(p,q) = d(G(a),G(b)) = \sum_{i=0}^0 d(G(t_i),G(t_{i+1}))
\]
where $a = t_0 < t_1 = b$. Assuming that the claim is true for $n$, let us prove it for $n+1$. Let $p,q \in S$, $C \subset S$ be any curve from $p$ to $q$ and $a = t_0 < t_1 < \cdots < t_{n+1} < t_{n+2} = b$. Since the image by $\pi_p$ of the minor arc from $p$ to $q$ is the segment from $p$ to $\pi_p q$ there exists a minimum $s$ such that $\pi_p C(s) = \pi_p G(t_{n+1})$. And, since the image by $\pi_q$ of the minor arc from $p$ to $q$ is the segment from $q$ to $\pi_q p$ there exists a maximum $t$ such that $\pi_q C(t) = \pi_q G(t_{n+1})$. It follows that $s \leq t$. Now we can write $C = C_s \cup C_{st} \cup C_t$, where $C_s = C_{as}$ and $C_t = C_{tb}$. We can rotate the curve $C_s$ through the axis defined by $p$ obtaining a curve $\widehat{C}_s$ from $p$ to $G(t_{n+1})$ and we can rotate the curve $C_t$ through the axis defined by $q$ obtaining a curve $\widehat{C}_t$ from $G(t_{n+1})$ to $q$. By the induction hypothesis, we have that
\[
\ell(\widehat{C}_s) \geq \sum_{i=0}^n d(G(t_i),G(t_{i+1}))
\]
Since
\[
\ell(\widehat{C}_t) \geq d(G(t_{n+1}),q) = d(G(t_{n+1}),G(t_{n+2}))
\]
it follows that
\[
\ell(C) = \ell(C_s) + \ell(C_{st}) + \ell(C_t)  \geq \ell(\widehat{C}_s) + \ell(\widehat{C}_t) \geq \sum_{i=0}^{n+1} d(G(t_i),G(t_{i+1}))
\]
since $\ell(C_{st}) \geq 0$. Now we fix $p,q \in S$. For any curve $C \subset S$ from $p$ to $q$ and any finite partition $a = t_0 < t_1 < \cdots < t_n < t_{n+1} = b$ of $[a,b]$, we have that
\[
\ell(C) \geq \sum_{i=0}^n d(G(t_i),G(t_{i+1}))
\]
Taking the supremum, it follows that $\ell(C) \geq \ell(G)$, showing that $G$ is a geodesic from $p$ to $q$.
\end{proof}

\end{document}